\documentclass{article}

\usepackage{amsmath,amssymb}
\usepackage{graphicx}
\usepackage{color}

\newtheorem{theorem}{Theorem}[section]%
\newtheorem{corollary}[theorem]{Corollary}%
\newtheorem{proposition}[theorem]{Proposition}%

\newtheorem{definition}[theorem]{Definition}%
\newtheorem{remark}[theorem]{Remark}%
\newtheorem{example}[theorem]{Example}%

\newcommand{\done}{\hfill $\Box$ }

\newcommand{\ls}[1]
    {\dimen0=\fontdimen6\the\font\lineskip=#1\dimen0
     \advance\lineskip.5\fontdimen5\the\font
     \advance\lineskip-\dimen0
     \lineskiplimit=0.9\lineskip
     \baselineskip=\lineskip
     \advance\baselineskip\dimen0
     \normallineskip\lineskip\normallineskiplimit\lineskiplimit
     \normalbaselineskip\baselineskip
     \ignorespaces}

\begin{document}

\bibliographystyle{abbrv}

\title{A counterexample shows that not every locally $L^0$--convex topology is necessarily induced by a family of $L^0$--seminorms}
\author{Mingzhi Wu$^\ast$ \quad Tiexin Guo\\
School of Mathematics and Statistics \\
Central South University\\
Changsha {\rm 410083}, China \\
Email: wumingzhi@csu.edu.cn, tiexinguo@csu.edu.cn\\
}

\date{}
 \maketitle

\thispagestyle{plain}
\setcounter{page}{1}

\begin{abstract}
 This paper constructs a counterexample showing that not every locally $L^0$--convex topology is necessarily induced by a family of $L^0$--seminorms. Random convex analysis is the analytic foundation for $L^0$--convex conditional risk measures, this counterexample, however, shows that a locally $L^0$--convex module is not a proper framework for random convex analysis. Further, this paper also gives a necessary and sufficient condition for a locally $L^0$--convex topology to be induced by a family of $L^0$--seminorms. Finally, we give some comments showing that based on random locally convex modules, we can establish a perfect random convex analysis to meet the needs of the study of $L^0$--convex conditional risk measures.


\end{abstract}

\ls{1.5}
\section{Introduction}

It is well known that classical convex analysis (see \cite{Ekeland}) is the analytic foundation for convex risk measures, cf.\cite{ADEH,Delbaen,FS}. However, classical convex analysis is not well suited to the study of conditional convex (or, $L^0$--convex) conditional risk measures defined on the spaces of unbounded financial positions. It is to overcome this obstacle that Filipovi\'c, Kupper and Vogelpoth \cite{FKV} presented the module approach to conditional risk. The key point in this module approach is to establish random convex analysis as the analytic foundation for $L^0$--convex conditional risk measures. To this, they introduced the notion of locally $L^0$--convex modules, as a module analogue of locally convex spaces. In the theory of locally convex spaces, it is a basic fact that every locally convex topology can be induced by a family of seminorms, it is to establish the module analogue of the basic fact that Filipovi\'c, Kupper and Vogelpoth \cite{FKV} further proved that the locally $L^0$--convex topology of every locally $L^0$--convex module can be induced by a family of $L^0$--seminorms, namely, Theorem 2.4 of \cite{FKV} as the basis for their whole paper. Unfortunately, this paper provides a counterexample to show that Theorem 2.4 of \cite{FKV} is wrong. Besides, in particular this paper gives a necessary and sufficient condition for a locally $L^0$--convex topology to be induced by a family of $L^0$--seminorms. Finally, we give some comments showing that based on random locally convex modules, we can establish a perfect random convex analysis to meet the needs of the study of $L^0$--convex conditional risk measures.

The remainder of this paper is organized as follows: in Section 2 we first recall some necessary terminology and notation; in Section 3 we construct the counterexample mentioned above and further give a necessary and sufficient condition for a locally $L^0$--convex topology to be induced by a family of $L^0$--seminorms; finally, in Section 4 we give some comments showing that based on random locally convex modules, we can establish a perfect random convex analysis to meet the needs of the study of $L^0$--convex conditional risk measures.

\section{Terminology and notation}

Let $(\Omega,{\mathcal F},P)$ be a probability space, $K$ the scalar field $\mathbb{R}$ of real numbers or $\mathbb{C}$ of complex numbers and $L^{0}(\mathcal{F},K)$ be the algebra of all equivalence classes of $K$--valued ${\mathcal F}$--measurable random variables on $\Omega$. Specially, $L^0=L^{0}(\mathcal{F},\mathbb{R})$. As usual, $L^{0}$ is partially ordered by $\xi\leqslant\eta$ iff $\xi^{0}(\omega)\leq\eta^{0}(\omega)$ for $P$--almost all $\omega\in \Omega$, where $\xi^0$ and $\eta^0$ are arbitrarily chosen representatives of $\xi$ and $\eta$, respectively. According to \cite{Dunford}, $(L^0,\leqslant)$ is a conditionally complete lattice. For a subset $A$ of $L^0$ with an upper bound (a lower bound), $\vee A$ (accordingly, $\wedge A$) stands for the supremum (accordingly, infimum) of $A$.
Let $\xi$ and $\eta$ be in $L^0$, we use $``\xi<\eta $ (or $\xi\leq\eta$) on $A"$ for $``\xi^{0}(\omega)<\eta^{0}(\omega)$ (resp., $\xi^{0}(\omega)\leq\eta^{0}(\omega)$) for $P$--almost all $\omega\in A"$, where $A\in\mathcal{F}$, $\xi^0$ and $\eta^0$ are a representative of $\xi$ and $\eta$, respectively.

Denote $L^{0}_{+}=\{\xi\in L^{0}\,|\,\xi\geqslant 0\}$ and $L^{0}_{++}=\{\xi\in L^{0}\,|\,\xi>0 \text{~on~$\Omega$}\}$.

${\tilde I}_A$ always denotes the equivalence class of $I_A$, where $A\in {\mathcal F}$ and $I_A$ is the characteristic function of $A$. For any $\xi\in L^{0}(\mathcal{F},K)$, $|\xi|$ denotes the equivalence class of $|\xi^0|: \Omega\to \mathbb{R}_+$ defined by $|\xi^0|(\omega)=|\xi^0(\omega)|$, where $\xi^0$ is an arbitrarily chosen representative of $\xi$.

For any $\varepsilon\in L^0_{++}$, denote $B_{\varepsilon}=\{\xi\in L^{0}(\mathcal{F},K)~|~|\xi|\leqslant\varepsilon\}$. Let $${\mathcal T}_c=\{V\subset L^{0}(\mathcal{F},K)~|~\text{for every $y\in V$ there exists $\varepsilon\in L^0_{++}$ such that $y+B_{\varepsilon}\subset V$}\},$$ then ${\mathcal T}_c$ is a Hausdorff topology on $L^{0}(\mathcal{F},K)$ such that $(L^{0}(\mathcal{F},K), {\mathcal T}_c)$ is a topological ring, namely the addition and multiplication operations are jointly continuous. D. Filipovi\'{c}, M. Kupper and N. Vogelpoth first introduced in \cite{FKV} this kind of topology and further pointed out that ${\mathcal T}_c$ is not necessarily a linear topology since the scalar multiplication mapping: $\alpha\mapsto \alpha x$ ($x$ is fixed) is no longer continuous in general. These observations led them to the study of a class of topological modules over the topological ring $(L^{0}(\mathcal{F},K), {\mathcal T}_c)$.

\begin{definition}(see \cite{FKV}).
A topological $L^{0}(\mathcal{F},K)$--module $(E, {\mathcal T})$ is an $L^{0}(\mathcal{F},K)$--module $E$ endowed with a topology ${\mathcal T}$ such that the addition and module multiplication operations:

\noindent~(i) $(E, {\mathcal T})\times (E, {\mathcal T})\to (E, {\mathcal T}),~(x_1,x_2)\mapsto x_1+x_2$ and\\
\noindent(ii) $(L^{0}(\mathcal{F},K),{\mathcal T_c})\times (E, {\mathcal T})\to (E, {\mathcal T}),~(\xi,x)\mapsto \xi x$

\noindent are continuous w.r.t. the corresponding product topologies.

\end{definition}

Locally $L^0$--convex topologies are defined as follows.

\begin{definition}(see \cite{FKV}).
For a topological $L^{0}(\mathcal{F},K)$--module $(E, {\mathcal T})$, the topology ${\mathcal T}$ is said to be locally $L^0$--convex if there is a neighborhood base ${\mathcal U}$ of $0\in E$ for which every $U\in {\mathcal U}$ is:

\noindent(i)  $L^0$--convex: $\xi x_1+(1-\xi)x_2\in U $ for all $x_1, x_2\in U$ and $\xi\in L^0$ with $0\leqslant \xi\leqslant 1$,\\
\noindent(ii) $L^0$--absorbent: for all $x\in E$ there is $\xi\in L^0_{++}$ such that $x\in \xi U$,\\
\noindent(iii) $L^0$--balanced: $\xi x\in U$ for all $x\in U$ and $\xi\in L^{0}(\mathcal{F},K)$ with $|\xi|\leqslant 1.$

\noindent In this case, $(E, {\mathcal T})$ is called a locally $L^0$--convex module.
\end{definition}

Given an $L^{0}(\mathcal{F},K)$--module $E$, an easy way to construct a locally $L^0$--convex topology on $E$ is by a family of $L^0$--seminorms on $E$. The notions of $L^0$--norms, $L^0$--seminorms and random locally convex modules were introduced by Guo before 2009 and random normed modules and random locally convex modules have been deeply developed under the $(\varepsilon, \lambda)$--topology, see Section 4 of this paper for the related terminology. Let us first recall the notion of $L^0$--seminorms.

\begin{definition}
Let $E$ be an $L^{0}(\mathcal{F},K)$--module, a function $\|\cdot\|: E\to L^0_+$ is called an $L^0$--seminorm on $E$ if:

\noindent(i)  $\|\xi x\|=|\xi|\|x\|$ for all $\xi\in L^{0}(\mathcal{F},K)$ and $x\in E$,\\
\noindent(ii) $\|x_1+x_2\|\leqslant \|x_1\|+\|x_2\|$ for all $x_1, x_2\in E$.

\noindent Furthermore, an $L^0$--seminorm $\|\cdot\|$ on $E$ is called an $L^0$--norm if $\|x\|=0$ implies $x=0$.
\end{definition}

\begin{proposition}\label{inducedtopology}(see \cite{FKV}).
Let $E$ be an $L^{0}(\mathcal{F},K)$--module and ${\mathcal P}$ a family of $L^0$-seminorms on $E$, for finite ${\mathcal Q}\subset {\mathcal P}$ and $\varepsilon\in L^0_{++}$ we define
$$U_{{\mathcal Q},\,\varepsilon}=\left\{x\in E~\right|\|x\|\leqslant \varepsilon,\forall\,\|\cdot\|\in {\mathcal Q}\},$$ then
$${\mathcal U}_{\mathcal P}=\{U_{{\mathcal Q},\,\varepsilon}~|~{\mathcal Q}\subset {\mathcal P}\text{~finite and~} \varepsilon\in L^0_{++}\}$$ forms a neighborhood base of $0$, of some locally $L^0$--convex topology on $E$, called the topology induced by ${\mathcal P}$.
\end{proposition}

\begin{definition}(see \cite{FKV}).
Let $E$ be an $L^{0}(\mathcal{F},K)$--module, the random gauge function $p_U:E\to L^0_+$ of an $L^0$--absorbent set $U\subset E$ is defined by $$p_U(x)=\wedge\{\xi\in L^0_+~|~x\in \xi U\},\quad\forall x\in E.$$
\end{definition}

It is proved in \cite{FKV} that: if $U\subset E$ is $L^0$--convex, $L^0$--absorbent and $L^0$--balanced, then the random gauge function $p_U$ is an $L^0$--seminorm on $E$ and
$p_U(x)=\wedge\{\xi\in L^0_{++}~|~x\in \xi U\}, \forall x\in E.$

For the subsequent use, we give two simple facts about random gauge function as follows.

\begin{proposition}\label{simplyfact}
Let $E$ be an $L^{0}(\mathcal{F},K)$--module, then we have the following:

\noindent(1).~For any $L^0$--seminorm $p$ on $E$, let $V=\{x\in E~|~p(x)\leqslant 1\}$, then $p_V=p$;

\noindent(2).~For any finite family ${\mathcal P}$ of $L^0$--seminorms on $E$ and $\varepsilon \in L^0_{++}$, let $U=\{x\in E~|~p(x)\leqslant \varepsilon,~\forall p\in {\mathcal P}\}$, then $\{x\in E~|~p_U(x)\leqslant 1\}=U$.
\end{proposition}
{\em Proof.} (1). For any given $x\in E$, we have $p_V(x)=\wedge\{\xi\in L^0_{++}~|~x\in \xi V\}=\wedge\{\xi\in L^0_{++}~|~\xi^{-1}x\in V \}=\wedge\{\xi\in L^0_{++}~|~\xi^{-1}p(x)\leqslant 1 \}=\wedge\{\xi\in L^0_{++}~|~p(x)\leqslant \xi \}=p(x)$.

(2). The inclusion $U\subset \{x\in E~|~p_U(x)\leqslant 1\}$ is clear from the definition, so it only needs to show the reverse inclusion. Let $x$ be an element of $E$ such that $p_U(x)\leqslant 1$, then, for any $p\in {\mathcal P}$ and $\delta\in L^0_{++}$ such that $x\in \delta U$, we have that $p(x)\leqslant \delta\cdot\vee\{p(y):y\in U\}\leqslant\delta\varepsilon$, therefore, $p(x)\leqslant \varepsilon\cdot\wedge\{\delta\in L^0_{++}~|~x\in \delta U\}=\varepsilon p_U(x)\leqslant\varepsilon$.
\hfill\done

\section{A counterexample and a necessary and sufficient condition for a locally $L^0$--convex topology to be induced by a family of $L^0$--seminorms}
 \label{}

 To construct an example of a locally $L^0$--convex topology which can not be induced by any family of $L^0$--seminorms, it is clear that we first need to find a new method (namely, not by use of $L^0$--seminorms as in Proposition \ref{inducedtopology}) to construct a locally $L^0$--convex topology. The following proposition is the basis for our method.

 \begin{proposition}\label{localbase}
 Let $E$ be an $L^{0}(\mathcal{F},K)$-module, $\mathcal{U}$ a family of $L^0$--convex, $L^0$--absorbent and $L^0$--balanced subsets of $E$ which satisfies the following three conditions:\\
 \noindent(1). For any $U_1, U_2\in \mathcal{U}$, there exists $U_3\in \mathcal{U}$ such that $U_3\subset U_1\cap U_2$,\\
 \noindent(2). For any $U\in \mathcal{U}$, there exists $V\in \mathcal{U}$ such that $V+V\subset U$,\\
 \noindent(3). For any $U\in \mathcal{U}$ and $\varepsilon\in L^0_{++}$, there exists $V\in \mathcal{U}$ such that $\varepsilon V\subset U$.\\
Let ${\mathcal T}=\{V\subset E~|~\text{for any $y\in V$, there is $U\in {\mathcal U} $ such that $y+U\subset V$}\}$, then $\mathcal T$ is a locally $L^0$--convex topology on $E$ and $\mathcal{U}$ is a neighborhood base of $\mathcal T$ at $0$.
 \end{proposition}
{\em Proof.}  It is easily seen that $\mathcal T$ is a topology on $E$. It remains to show that $(E,{\mathcal T})$ is a topological $L^{0}(\mathcal{F},K)$--module.

For any $x, y\in E$ and $U\in {\mathcal U}$, by (2), there exists $V\in {\mathcal U}$ such that $V+V\subset U$, it follows that $(x+V)+(y+V)=(x+y)+(V+V)\subset (x+y)+U$, thus the addition operation is continuous.

For any $t\in L^{0}(\mathcal{F},K),\,x\in E$ and $U\in {\mathcal U}$, by (2), there exists $V\in {\mathcal U}$ such that $V_1+V_1\subset U$. Since $V_1$ is $L^0$-absorbent, there exists $\varepsilon\in L^0_{++}$ such that $\varepsilon x\in V_1$. Let $\delta=\varepsilon +|t|$, according to (3), there exists $V_2\in {\mathcal U}$ such that $\delta V_2\subset V_1$. Noting that $V_1$ is $L^0$--balanced, we can further obtain that $B_\varepsilon x:=\{\xi x:~|\xi|\leqslant \varepsilon\}\subset V_1$ and $(t+B_\varepsilon)V_2\subset B_\delta V_2\subset V_1$, therefore, $$(t+B_\varepsilon)(x+V_2)=tx+B_\varepsilon x+(t+B_\varepsilon)V_2\subset tx+V_1+V_1\subset tx+U,$$ which means that the module multiplication operation is continuous.
\hfill\done

Now we can give an example of a locally $L^0$--convex module $(E, {\mathcal T})$ for which the topology ${\mathcal T}$ can not be induced by any family of $L^0$--seminorms on $E$.

For the sake of convenience, we first recall a notation here. For a subset $C$ of an $L^0$--module $E$, we denote by
 $$ span_{L^0}(C):=\left\{\sum^n_{i=1}x_ic_i~\left|~c_i\in C, x_i\in L^0, 1\leq i\leq n, n\in \mathbb{N}\right.\right\}$$
 the $L^0$--submodule of $E$ generated by $C$.

\begin{example}\label{counterexample}
 Let $\Omega$ be the set of all positive integers, ${\mathcal F}$ the $\sigma$--algebra of all the subsets of $\Omega$, the probability $P$ on $(\Omega, {\mathcal F})$ defined by $P(\{j\})=2^{-j},\,\forall j\in \Omega$. Let $E=L^0$, $M=span_{L^0}\{{\tilde I}_{\{j\}}:j\in \Omega\}$, and for any $\varepsilon\in L^0_{++}$, recall that $B_{\varepsilon}=\{x\in E~|~|x|\leqslant \varepsilon\}$, then both $M$ and $B_{\varepsilon}$ are $L^0$--convex and $L^0$--balanced. In addition, $B_{\varepsilon}$ is $L^0$-absorbent. Let $U_{\varepsilon}=M+B_{\varepsilon}$, then $U_{\varepsilon}$ is an $L^0$--convex, $L^0$-absorbent and $L^0$--balanced subset of $E$. Moreover, the family ${\mathcal U}=\{U_\varepsilon: \varepsilon\in L^0_{++}\}$ satisfies the three conditions in Proposition \ref{localbase}, in fact, for any $\varepsilon,\delta\in L^0_{++}$, we have $U_{\varepsilon\wedge \delta}\subset U_{\varepsilon}\cap U_{\delta}$, $U_{\varepsilon/2}+U_{\varepsilon/2}\subset U_{\varepsilon}$ and $\varepsilon U_{\delta/\varepsilon}\subset U_{\delta}$. Therefore, ${\mathcal U}$ forms a neighborhood base of $0$ of a locally $L^0$--convex topology (denoted by ${\mathcal T}$) on $E$. We have that the topology ${\mathcal T}$ can not be induced by any family of $L^0$--seminorms on $E$.
\end{example}
{\em Proof.}
Fix one $\varepsilon\in L^0_{++}$, we calculate the random gauge function $p_{U_\varepsilon}$. For each $x\in E$ and each $j\in \Omega$, since $${\tilde I}_{\{j\}}x\in span_{L^0}\{{\tilde I}_{\{j\}}\}\subset M=\delta M\subset \delta U_\varepsilon,\,\forall \delta\in L^0_{++},$$ we have that ${\tilde I}_{\{j\}}p_{U_\varepsilon}(x)=p_{U_\varepsilon}({\tilde I}_{\{j\}}x)\leqslant \wedge L^0_{++}=0$, namely, $p_{U_\varepsilon}(x)=0,\forall x\in E$.

 We show that ${\mathcal T}$ can not be induced by any family of $L^0$--seminorms by contradiction. Assume that ${\mathcal P}$ is a family of $L^0$--seminorms on $E$ which induces the topology ${\mathcal T}$, then for any $p\in {\mathcal P}$, there exists $U_\varepsilon$ such that $U_\varepsilon\subset V:=\{x\in E:~p(x)\leqslant 1\}$, according to Proposition \ref{simplyfact}, we have that $p=p_V\leqslant p_{U_\varepsilon}=0$, thus $p$ must also be zero. It follows that the family of $L^0$--seminorms ${\mathcal P}$ is actually a singleton $\{0\}$. Hence, the locally $L^0$--convex topology induced by ${\mathcal P}$ is the trivial chaos topology which consists of $\emptyset$ and $E$. However, we claim that $M$ is a proper ${\mathcal T}$--closed $L^0$--submodule of $E$, which implies that ${\mathcal T}$ is not trivial. This is a contradiction.

 It remains to show the claim to complete the proof. Clearly, $M$ is a proper subset of $E$, we only need to show that $M$ is ${\mathcal T}$--closed. Since the ${\mathcal T}$--closure of $M$ equals $\bigcap\{M+U_\varepsilon: \varepsilon\in L^0_{++}\}=\bigcap\{M+M+B_\varepsilon: \varepsilon\in L^0_{++}\}=\bigcap\{M+B_\varepsilon: \varepsilon\in L^0_{++}\}$, we need to show that $\bigcap\{M+B_\varepsilon: \varepsilon\in L^0_{++}\}=M$, that is to say, for an arbitrarily given $x\in E$ which is not in $M$, we need to find an $\varepsilon\in L^0_{++}$ such that $x$ is not in $M+B_\varepsilon$. To this end, let $S=\{j\in \Omega:x(j)\neq 0\}$, then $S$ is an infinite set, otherwise, if $S$ is finite, then $x\in span_{L^0}\{{\tilde I}_{\{j\}}:j\in S\}\subset M$. Define $\varepsilon\in L^0_{++}$ by
$$\varepsilon(j)=\left\{
                  \begin{array}{ll}
                   {{1}\over{2}}|x(j)|, & \hbox{$j\in S$;} \\
                    1, & \hbox{$j\in \Omega\setminus S$,}
                  \end{array}
                \right.$$
then $x$ is not in $M+B_\varepsilon$. In fact, if $x=m+y$ for some $m\in M$ and $y\in B_\varepsilon$, then for each $j\in S$, we have that $|m(j)|=|x(j)-y(j)|\geq |x(j)|-|y(j)|\geq |x(j)|-\varepsilon(j)={{1}\over{2}}|x(j)|>0$. However, according to the definition of $M$, $\{j\in\Omega: m(j)\neq 0\}$ should be a finite set.
\hfill \done

 \begin{remark}
 In the above example, one can easily see that the topology ${\mathcal T}$ is not Hausdorff since the closure of $\{0\}$ equals $\bigcap\{0+U_\varepsilon: \varepsilon\in L^0_{++}\}=\bigcap\{M+B_\varepsilon: \varepsilon\in L^0_{++}\}=M$. By considering the quotient $L^0$--module $E/M$ and the quotient topology (denoted by ${\mathcal T}_\Pi$) on it, we get an example of a locally $L^0$--convex module $(E/M,{\mathcal T}_\Pi )$ which is Hausdorff but the topology cannot be induced by any family of $L^0$--seminorms.

\end{remark}

 \begin{remark}\label{remark}
In fact, for an arbitrarily given probability space $(\Omega,{\mathcal F},P)$ which is not essentially generated by finitely many $P$-atoms, by making a slight modification to Example \ref{counterexample}, we can always give a locally $L^0$--convex topology for $L^0$ which cannot induced by any family of $L^0$--seminorms. In fact, in such a situation, there exists a countable partition $\{A_n:n\in {\mathbb N}\}$ of $\Omega$ to ${\mathcal F}$ such that each $A_n$ has positive probability, we can verify this fact as follows: first, $(\Omega,{\mathcal F},P)$ has at most countably many disjoint $P$-atoms, which is denoted by $\{B_n: n=1,2,\dots,n_0\}$, where $n_0$ is a positive integer or $n_0=+\infty$, further, let $\Omega^\prime=\Omega\setminus {\cup_{n=1}^{n_0}B_n}$, then $P(\Omega^\prime)>0$ and $\Omega^\prime$ does not include any $P$-atoms, and hence there is a countable disjoint family $\{C_n\in {\mathcal F}: n\in \mathbb N\}$ such that $P(C_n)=\frac{P(\Omega^\prime)}{2^n}$ for each positive integer $n$, now let $\{A_n: n\in {\mathbb N}\}=\{B_n: n=1,2,\dots,n_0\}\cup \{C_n: n\in \mathbb N\}$, it is clear that $\{A_n: n\in {\mathbb N}\}$ is a countable partition of $\Omega$ to ${\mathcal F}$. We only need to set $M=span_{L^0}\{{\tilde I}_{A_n}:n\in {\mathbb N}\}$ and still let $U_\varepsilon=M+B_\varepsilon$ for each $\varepsilon\in L^0_{++}$, then the locally $L^0$--convex topology for $L^0$ induced by the local neighborhood base $\{U_\varepsilon:\varepsilon\in L^0_{++}\}$ meets our need, one can complete the verification only by replacing ${\tilde I}_{\{n\}}$ with ${\tilde I}_{A_n}$ in the proof of Example \ref{counterexample}. Finally, we know from the referee's report that J. M. Zapata has also, independently, obtained a similar result in \cite{Zapata}, where J. M. Zapata only assumes that $\Omega$ has a countable partition of positive probabilities, so he also obtained an enough general counterexample although \cite{Zapata} did not give the details of the verification of his counterexample and he did not consider a counterexample with the Hausdorff property.
 \end{remark}

Let $(E,{\mathcal T})$ be a locally $L^0$--convex module. Assume that ${\mathcal U}$ is a neighborhood base of $0\in E$ such that each $U\in {\mathcal U}$ is $L^0$--convex, $L^0$--absorbent and $L^0$--balanced. As pointed out in \cite{Guo-Zhao-Zeng}, if the inclusion relation
  \begin{equation*}
   \{x\in E~|~p_U(x)<1~\text{on~$\Omega$}\}\subset U\subset \{x\in E~|~p_U(x)\leqslant 1\}
  \end{equation*}
 holds true for each $U\in {\mathcal U}$, in particular if the relation
  \begin{equation}\label{equality}
   U=\{x\in E~|~p_U(x)\leqslant 1\}
  \end{equation}
 holds true for each $U\in {\mathcal U}$, then ${\mathcal T}$ can be induced by the family of $L^0$--seminorms $\{p_U:~U\in {\mathcal U}\}$. Conversely, assume that ${\mathcal P}$ is a family of $L^0$--seminorms on $E$ which induces ${\mathcal T}$, then for the neighborhood base ${\mathcal U}_{\mathcal P}$ as in Proposition \ref{inducedtopology}, according to Proposition \ref{simplyfact}, each $U\in {\mathcal U}_{\mathcal P}$ satisfies the relation \eqref{equality}. What is more important is that the relation \eqref{equality} helps us look for a necessary and sufficient condition for a locally $L^0$--convex topology to be induced by a family of $L^0$--seminorms. In fact, a sufficient condition was already given in \cite{Guo-Zhao-Zeng}, which requires that $U$ has the countable concatenation property.  Let us recall this important notion as follows.

 It should be pointed out that when Guo introduced the countable concatenation property for a subset of an $L^{0}(\mathcal{F},K)$--module $E$, $E$ is assumed to have the following property:\\

\noindent(C)\quad for any $x,~y\in E$, if there is a countable partition $\{A_n,n\in \mathbb{N}\}$ of $\Omega$ to ${\mathcal F}$ such that ${\tilde I}_{A_n}x={\tilde I}_{A_n}y$ for each $n\in N$, then $x=y$.\\

\noindent Guo proved in \cite{Guotx-relation} that every random locally convex module has the above property (C), however, up to now, no one has ever shown that every locally $L^0$--convex module necessarily has the property (C), in this paper for our purpose it is convenient to quit the requirement of the above property (C) and rephrase the notion of the countable concatenation property in a slightly wide sense as follows.

\begin{definition} Let $E$ be an $L^{0}(\mathcal{F},K)$--module. A sequence $\{x_n, n\in \mathbb{N}\}$ in $E$ is countably concatenated in $E$ with respect to a countable partition $\{A_n,n\in \mathbb{N} \}$ of $\Omega$ to ${\mathcal F}$ if there is $x\in E$ such that ${\tilde I}_{A_n}x={\tilde I}_{A_n}x_n$ for each $n\in \mathbb{N}$, in which case we denote the set of all such $x$ by $\sum^{\infty}_{n=1}{\tilde I}_{A_n}x_n$, namely, $$\sum^{\infty}_{n=1}{\tilde I}_{A_n}x_n=\{x\in E~|~\text{$x\in E$ such that ${\tilde I}_{A_n}x={\tilde I}_{A_n}x_n$ for each $n\in \mathbb{N}$}\}.$$ A subset $G$ of $E$ is said to have the countable concatenation property if each sequence $\{x_n, n\in \mathbb{N}\}$ in $G$ is countably concatenated in $E$ with respect to an arbitrary countable partition $\{A_n,n\in \mathbb{N}\}$ of $\Omega$ to $\mathcal{F}$ and $\sum^{\infty}_{n=1}{\tilde I}_{A_n}x_n\subset G$.
\end{definition}

In fact, Guo et al \cite{Guo-Zhao-Zeng} already proved Proposition \ref{prop} and Corollary \ref{corol} below.

\begin{proposition}\label{prop}(See \cite{Guo-Zhao-Zeng}).
Let $(E,{\mathcal T})$ be a locally $L^0$--convex module and $U$ an $L^0$--convex, $L^0$--absorbent and $L^0$--balanced subset with the countable concatenation property. Then $\{x\in E~|~p_U(x)<1~\text{on~$\Omega$}\}\subset U\subset \{x\in E~|~p_U(x)\leqslant 1\}$.
\end{proposition}

\begin{corollary}\label{corol}(See \cite{Guo-Zhao-Zeng}).
Let $(E,{\mathcal T})$ be a locally $L^0$--convex module. If there exists a neighborhood base ${\mathcal U}$ of $0\in E$ such that each $U\in {\mathcal U}$ is an $L^0$--convex, $L^0$--absorbent and $L^0$--balanced subset with the the countable concatenation property, then ${\mathcal T}$ can be induced by the family of $L^0$--seminorms $\{p_U:~U\in {\mathcal U}\}$.
\end{corollary}

In order to give a necessary and sufficient condition for a locally $L^0$--convex topology to be induced by a family of $L^0$--seminorms, we need Definition \ref{def} below, which is based on Guo's earlier work \cite{GuoarX}, where the notion of the relative countable concatenation property was first considered although the terminology  of the relative countable concatenation property did not occur in \cite{GuoarX}.

\begin{definition}\label{def}

A subset $G$ of an $L^{0}(\mathcal{F},K)$--module $E$ is said to have the relative countable concatenation property if for a sequence $\{x_n, n\in \mathbb{N}\}$ in $G$ and a countable partition $\{A_n,n\in \mathbb{N}\}$ of $\Omega$ to $\mathcal{F}$, we always have $\sum^{\infty}_{n=1}{\tilde I}_{A_n}x_n\subset G$ whenever $\{x_n, n\in \mathbb{N}\}$ is countably concatenated in $E$ with respect to $\{A_n,n\in \mathbb{N}\}$.
\end{definition}

\begin{remark}
For an arbitrary $L^{0}(\mathcal{F},K)$--module $E$, $E$ need not have the countable concatenation property, but it is clear that $E$ as a subset of itself always has the relative countable concatenation property. Furthermore, one can easily see that``$G$ has the relative countable concatenation property'' is the same as ``$G$ has the countable concatenation property'' for every subset $G$ when $E$ has the countable concatenation property. In Definition \ref{def}, the adjective ``relative'' means that whenever a sequence $\{x_n, n\in \mathbb{N}\}$ in $G$ is countably concatenated ``in $E$'' with respect to a countable partition $\{A_n,n\in \mathbb{N}\}$, then this sequence must be countably concatenated ``in $G$'' with respect to the countable partition $\{A_n,n\in \mathbb{N}\}$.
\end{remark}

\begin{theorem}
  Let $(E,{\mathcal T})$ be a topological $L^{0}(\mathcal{F},K)$--module, then the following two statements are equivalent to each other:\\
  \noindent (1). The topology ${\mathcal T}$ can be induced by a family of $L^0$--seminorms on $E$;\\
  \noindent (2). There exists a neighborhood base ${\mathcal U}$ of $0\in E$ such that each element $U\in {\mathcal U}$ is an $L^0$--convex, $L^0$--absorbent and $L^0$--balanced subset with the relative countable concatenation property.
\end{theorem}
{\em Proof.}
$(1)\Rightarrow(2)$. Assume that the topology ${\mathcal T}$ can be induced by a family of $L^0$--seminorms ${\mathcal P}$ on $E$. For any $\varepsilon\in L^0_{++}$ and finite ${\mathcal Q}\subset {\mathcal P}$, let $U_{{\mathcal Q},\,\varepsilon}=\left\{x\in E~\right|~\|x\|\leqslant \varepsilon,\forall\,\|\cdot\|\in {\mathcal Q}\}$, then ${\mathcal U}:=\{U_{{\mathcal Q},\,\varepsilon}~|~{\mathcal Q}\subset {\mathcal P}\text{~finite and~} \varepsilon\in L^0_{++}\}$ is a neighborhood base of $0$. We need only to show that each $U_{{\mathcal Q},\,\varepsilon}$ has the relatively countable concatenation property. To this end, assume that $\{x_n, n\in \mathbb{N}\}$ is a sequence in $U_{{\mathcal Q},\,\varepsilon}$ which is countably concatenated in $E$ with respect to a countable partition $\{A_n,n\in \mathbb{N}\}$ of $\Omega$ to $\mathcal{F}$. If $x\in \sum^{\infty}_{n=1}{\tilde I}_{A_n}x_n$, then for each $\|\cdot\|\in {\mathcal Q}$, we have that $$\|x\|=(\sum^{\infty}_{n=1}{\tilde I}_{A_n})\|x\|=\sum^{\infty}_{n=1}{\tilde I}_{A_n}\|x\|=\sum^{\infty}_{n=1}\|{\tilde I}_{A_n}x\|=\sum^{\infty}_{n=1}\|{\tilde I}_{A_n}x_n\|\leqslant \sum^{\infty}_{n=1}{\tilde I}_{A_n}\varepsilon=\varepsilon,$$
thus $x\in U_{{\mathcal Q},\,\varepsilon}$ and this in turn implies that $\sum^{\infty}_{n=1}{\tilde I}_{A_n}x_n\subset U_{{\mathcal Q},\,\varepsilon}$.

$(2)\Rightarrow (1)$. If (2) holds true, we will show that the topology ${\mathcal T}$ can be induced by the family of $L^0$--seminorms $\{p_U:~U\in {\mathcal U}\}$. It suffices to show that $\{x\in E~|~p_U(x)<1~\text{on ~$\Omega$}\}\subset U$. Let $x$ be a element of $E$ such that $p_U(x)<1~\text{on ~$\Omega$}$. Let ${\mathcal A}=\{A\in {\mathcal F}~|~{\tilde I}_Ax\in U\}$ and $A=ess.sup{\mathcal A}$, according to \cite[Proposition 2.25]{FKV}, $p_U(x)\geqslant 1 \text~{on~} A^c$, thus $P(A^c)=0$, namely, $P(A)=1$. Since $U$ is $L^0$--convex, ${\mathcal A}$ is directed upward and there exists an increasing sequence $\{A_n, n\in \mathbb{N}\}$ in ${\mathcal A}$ such that $A=\bigcup_{n\in \mathbb{N}}A_n$. Let $B_1=A_1\cup A^c, B_2=A_2\setminus A_1, \dots, B_n=A_n\setminus A_{n-1},\dots$, then $\{B_n,n\in \mathbb{N}\}$ is a countable partition of $\Omega$ to ${\mathcal F}$, and for each $n$ we have that ${\tilde I}_{B_n}x={\tilde I}_{B_n}({\tilde I}_{A_n}x)\in U$ since $U$ is $L^0$--balanced. Thus $x\in \sum^{\infty}_{n=1}{\tilde I}_{B_n}{\tilde I}_{B_n}x\subset U$ since $U$ has the relative countable concatenation property.\hfill \done

\begin{remark}
Recently, we know from the referee's report that J.M.Zapata \cite{Zapata} also independently presented the notion of being closed under the countable concatenation operation and further established the characterization theorem for a locally $L^0$--convex topology to be induced by a family of $L^0$-seminorms, see Page 6 and Theorem 2.1 in \cite{Zapata}. One can easily see that the relative countable concatenation property in our paper is the same as being closed under the countable concatenation operation, so J.M.Zapata's Theorem 2.1 is also the same as our Theorem 3.10.
\end{remark}

 \section{Concluding remarks}
Random convex analysis was first studied in \cite{FKV}, where a locally $L^0$--convex module is chosen as the space framework for random convex analysis. However, Example \ref{counterexample} and Remark \ref{remark} show that a locally $L^0$--convex module is not a proper space framework for random convex analysis. First, according to the notion of a locally $L^0$--convex module, only the locally $L^0$--convex topology can be employed, whereas for a random locally convex module $(E,\mathcal{P})$, the two kinds of topologies, namely, the $(\varepsilon, \lambda)$--topology and the locally $L^0$--convex topology can be used. Second, the locally $L^0$--convex topology seems more intuitive, but the $(\varepsilon, \lambda)$--topology is more natural from probability theory, in particular a random locally convex module absorbs both the advantages of this two kinds of topologies and the work of \cite{Guo-Zhao-Zeng} shows that a random locally convex module seems to be a more proper framework for random convex analysis. It is to overcome the disadvantage of a locally $L^0$--convex module that Guo et al \cite{Guo-Zhao-Zeng} choose a random locally convex module as a space framework to establish random convex analysis. Let us recall the notion of a random locally convex module as follows.

\begin{definition}($See$ \cite{TXG-Sur,TXG-strict}). An ordered pair $(E,\mathcal{P})$ is called a random locally convex module (briefly, an $RLC$ module) over $K$ with base $(\Omega,\mathcal{F},P)$ if $E$ is an $L^{0}(\mathcal{F},K)$--module and $\mathcal{P}$ a family of $L^0$--seminorms on $E$ such that $\vee\{\|x\|:\|\cdot\|\in\mathcal{P}\}=0$ iff $x=\theta$ (the null element of $E$).

\end{definition}

Let $\mathcal{P}$ be a family of $L^0$--seminorms on an $L^{0}(\mathcal{F},K)$--module $E$ and ${\mathcal P}_f$ the family of all finite subsets $\mathcal{Q}$ of ${\mathcal P}$. For each $\mathcal{Q}\in {\mathcal P}_f$, the $L^0$--seminorm $\|\cdot\|_\mathcal{Q}$ is defined by $\|x\|_\mathcal{Q}=\vee\{\|x\|:~\|\cdot\|\in \mathcal{Q}\},~\forall x\in E$.

\begin{definition}($See$ \cite{TXG-Sur,TXG-strict}.) Let $(E, {\mathcal P})$ be an $RLC$ module over $K$ with base $(\Omega, {\mathcal F}, P)$. For any positive numbers $\varepsilon$ and $\lambda$ with $0<\lambda<1$ and $\mathcal{Q}\in {\mathcal P}_f$, let $N_{\theta}(\mathcal{Q}, \varepsilon, \lambda)=\{x\in E~|~P\{\omega\in \Omega~|~\|x\|_\mathcal{Q}(\omega)<\varepsilon\}>1-\lambda\}$, then $\{N_{\theta}(\mathcal{Q}, \varepsilon, \lambda)~|~\mathcal{Q}\in {\mathcal P}_f, \varepsilon >0, 0<\lambda<1\}$ forms a local base at $\theta$ of some Hausdorff linear topology on $E$, called the $(\varepsilon, \lambda)$--topology induced by ${\mathcal P}$.
\end{definition}

 Since Guo's paper \cite{Guotx-relation}, random metric theory has come into such a model that random locally convex modules are developed by simultaneously considering the above two kinds of topologies, which makes random metric theory deeply developed (see, e.g. \cite{Guotx-recent,Zhao-Guo,Guo-shi,TXG-YJY,Wu1,Wu2,Guotx-onsome}) since the connection between basic results derived from the two kinds of topologies has been established in \cite{Guotx-relation}. Based on these deep advances, a complete random convex analysis has been developed in \cite{Guo-Zhao-Zeng} and some concrete applications of random convex analysis to $L^0$--convex conditional risk measures are also given in \cite{GZZ}.

 \vspace{0.8cm}

{\noindent\bf Acknowledgements.}

\vspace{0.5cm}

The first author is supported by Central South University Postdoctoral Science Foundation. The second author is supported by National Natural Science Foundation of China (Grant No. 11171015). The authors would like to thank the referee for pointing out to us the reference \cite{Zapata}, which makes us, for the first time, see J.M.Zapata's work \cite{Zapata}. Although \cite{Zapata} appeared on arXiv on April 29, 2014, sooner than the submitting time of our paper, one of the authors, namely, Prof.Guo told J.M.Zapata (on April 2, 2014) that we had obtained Counterexample 3.2 of our paper in an email to J.M.Zapata which did not include any details of our paper, in fact we then had obtained all the results of our paper but Remark \ref{remark}(Remark \ref{remark} is added after reading \cite{Zapata}). Thus we can say that our results are independent of \cite{Zapata}.


\end{document}